%% file: arXiv.tex
\documentclass[journal]{IEEEtran}

\usepackage{cite}
\usepackage{amssymb,amsmath, amsfonts, amsthm}
\usepackage[caption=false,font=footnotesize]{subfig}
\usepackage{tikz,pgfplots}
\usetikzlibrary{shapes.geometric, decorations.pathmorphing, positioning, arrows, arrows.meta, patterns, shapes, calc, fit, backgrounds, external}

\usepackage{graphicx}
\graphicspath{{./figures/}}
\DeclareGraphicsExtensions{.pdf,.jpeg,.png,.eps}
\usepackage[dvipsnames]{xcolor}

\input{lvc_style}

\begin{document}

\title{Multitrace M\"{u}ller Boundary Integral Equation for Electromagnetic Scattering by Composite Objects}

\author{Van~Chien~Le,~\IEEEmembership{Member,~IEEE,} and~Kristof~Cools,~\IEEEmembership{Member,~IEEE}
\thanks{Manuscript received April 19, 2005; revised August 26, 2008; accepted 19 May 2014. Date of publication 12 June 2014; date of current version 9 July 2014. This work was supported by the European Research Council (ERC) under the European Union's Horizon 2020 Research and Innovation Program under Grant 101001847. \textit{(Corresponding author: Van Chien Le.)}}%
\thanks{The authors are with IDLab, Department of Information Technology at Ghent University -- imec, 9000 Ghent, Belgium (e-mail: vanchien.le@ugent.be, kristof.cools@ugent.be).}
}


\maketitle
 
\begin{abstract}
     This paper introduces a boundary integral equation for time-harmonic electromagnetic scattering by composite dielectric objects. The formulation extends the classical M\"{u}ller equation to composite structures through the global multi-trace method. The key ingredient enabling this extension is the use of the Stratton--Chu representation in complementary region, also known as the extinction property, which augments the off-diagonal blocks of the interior representation operator. The resulting block system is composed entirely of second-kind operators. A Petrov--Galerkin (mixed) discretization using Rao--Wilton--Glisson trial functions and Buffa--Christiansen test functions is employed, yielding linear systems that remain well conditioned on dense meshes and at low frequencies without the need for additional stabilization. This reduces computational costs associated with matrix-vector multiplications and iterative solving. Numerical experiments demonstrate the accuracy of the method in computing field traces and derived quantities.
\end{abstract}



\section{Introduction}

\IEEEPARstart{T}{he} accurate modelling of composite electromagnetic systems is essential for the analysis and design of modern devices in antennas, propagation, microelectronics, and photonics. However, computing electromagnetic fields in the presence of heterogeneous materials and geometrically complex, multi-scale structures remains challenging. Among several numerical methods, boundary integral equation (BIE) formulations have emerged as a powerful and widely used approach. Their main advantage lies in the reduction of problem dimensionality, since only the interfaces of the domain require discretization, leading to a significant decrease in the number of unknowns compared to volumetric-based methods. Moreover, BIEs inherently satisfy the radiation condition, making them particularly effective for open-region scattering and radiation problems.

BIE formulations for composite electromagnetic problems differ primarily in how they enforce transmission conditions at interfaces between adjacent materials. Many approaches have been introduced over the years, most of them can be classified into the single-trace, local multi-trace, and global multi-trace formulations. In single-trace formulations, a single pair of unknowns representing the tangential traces of the electric and magnetic fields is sought on each interface, thereby enforcing field continuity in a direct manner. In contrast, multi-trace methods introduce two independent pairs of traces on each interface, one for each neighbouring subdomain. Local multi-trace formulations combine the self-representing property of the tangential traces in each domain with the continuity conditions, resulting in a uniquely solvable system for the traces. Global multi-trace formulations, on the other hand, impose these conditions weakly by introducing a conceptual gap between adjacent regions. These distinctions lead to formulations defined on different trace spaces, and the differences become particularly pronounced in the presence of complex junctions where three or more regions meet.

Over the past decades, considerable efforts have been devoted to developing efficient and robust BIE solvers for composite electromagnetic scattering problems. First-kind formulations, despite their ill-conditioning, remain widely used due to their high accuracy and natural extension to composite structures. Single-trace extensions of the classical Poggio--Miller--Chang--Harrington--Wu--Tsai (PMCHWT) equation have been investigated in \cite{CJS+2014,KBH+2019}, although stabilizing these formulations on dense meshes remains a significant challenge, especially in the presence of junctions. A quasi-local Calder\'{o}n preconditioner recently proposed in \cite{Cools2026} partially mitigates this difficulty and achieves nearly mesh-independent Krylov convergence. Local multi-trace methods have also been examined in \cite{PLL2013,Peng2015,ZHH+2018}, but their efficiency depends strongly on the discretization and on how the transmission conditions are enforced. The global multi-trace PMCHWT formulation was independently explored in \cite{CCZ+2003,YKJ2011,CH2012}, where its well-posedness was established, and low-frequency stabilization and dense-mesh preconditioning strategies were introduced. This framework was later extended to structures involving screens and perfect electric conductors in \cite{LBG+2022}. Nevertheless, a first-kind formulation that is simultaneously efficient, accurate, and robust for general composite configurations across broad frequency ranges is still lacking.

In contrast to first-kind BIEs, second-kind formulations for dielectric scattering, such as the M\"{u}ller equation \cite{Muller1969} and the combined N formulation (CNF) \cite{YTJ2005}, are intrinsically well-conditioned, leading to fast convergence of iterative solvers, even on dense meshes. Among them, the M\"{u}ller equation is generally preferred because it remains well-conditioned in the low-frequency regime, although the CNF may offer slightly better performance at certain moderate frequencies and material contrasts \cite{YTJ2005}. The superior dense-mesh and low-frequency conditioning of the M\"{u}ller equation stems from the cancellation of hyper-singular contributions, which regularizes the integral kernels. This regularity also permits a broad range of discretization schemes, including those based on non-conforming and discontinuous basis functions \cite{YT2005,UTR2011,CHS2017}. Despite these advantages, the M\"{u}ller formulation, like other second-kind BIEs, was historically regarded as significantly less accurate than its first-kind counterparts \cite{YTJ2005,YTJ2008,UTR2011}. This perception changed with the introduction of dual basis functions, such as the Buffa--Christiansen (BC) \cite{BC2007} and Chen--Wilton (CW) \cite{CW1990} functions, which enabled conforming mixed discretizations that preserve the favorable conditioning of the M\"{u}ller formulation while substantially improving both accuracy and low-frequency behavior, in many cases to a level comparable to that of the first-kind PMCHWT equation \cite{YJN2011,YJN2013,BCA+2014,YKJ2016}. Furthermore, these mixed discretizations belong to the Petrov--Galerkin framework and hence, when uniformly stable, they achieve optimal convergence with respect to spatial discretization, as established in \cite{LC2024b}.

Unfortunately, in composite settings involving junctions, it remained unclear how to cancel the hyper-singularities associated with all participating single-layer operators \cite{YTJ2005,TY2006,SLY+2009}. The single-trace M\"{u}ller formulation introduced in \cite{CHS2017} represents a significant step forward by enabling this singularity cancellation through the incorporation of the multi-potential representation. Its discontinuous Galerkin discretization using piecewise-constant basis functions yields rapidly convergent systems. Nevertheless, the non-conforming nature of this scheme introduces low-frequency instabilities and degrades the accuracy of secondary quantities, such as far-field evaluations. More broadly, this issue reflects a fundamental limitation of single-trace formulations for composite structures: no discrete dual of single-trace energy spaces exists on interfaces that contain junctions. As a result, constructing effective Calder\'{o}n preconditioners for first-kind BIEs and conforming mixed discretizations of second-kind equations within the single-trace framework becomes extremely challenging.

In this paper, we introduce a global multi-trace M\"{u}ller formulation for electromagnetic scattering by arbitrary composite dielectric objects. In contrast to the single-trace approach in \cite{CHS2017}, the proposed formulation is derived through the use of the Stratton--Chu representation in the complementary region, i.e., exploiting the extinction property. The off-diagonal block $(j, k)$ of the interior boundary representation operator is augmented with the tangential traces, taken on the boundary of the $j$-th subdomain, of the potential operators associated with the $k$-th subdomain. This augmentation enables a consistent linear combination of the exterior and interior representation formulas that eliminates the hyper-singularities in all block-matrix entries, resulting in a system composed entirely of second-kind operators. A key advantage of the proposed formulation is that it admits a conforming mixed discretization, based on trial and dual basis functions defined on the full boundary of each subdomain. This provides, for the first time, a second-kind multi-trace formulation for composite dielectric scattering that is both conforming and applicable in the presence of junctions, thereby overcoming the instabilities inherent in existing first-kind or non-conforming second-kind schemes. Numerical experiments demonstrate that the new formulation achieves accurate field traces and derived quantities, while exhibiting significantly improved conditioning at low frequencies and on dense meshes. Moreover, despite the introduction of additional boundary integral operators (BIOs) and dual test functions, the overall computational cost remains competitive with that of first-kind formulations, including Calder\'{o}n-preconditioned (CP) variants, when the number of scattering components is moderate.

An important motivation for developing the global multi-trace M\"{u}ller equation lies in its potential extension to time-domain (TD) problems. In this regime, PMCHWT formulations are known to suffer from late-time instabilities \cite{LGC+2025}. Although stabilization strategies for TD-PMCHWT equations have been proposed, they typically lead to substantially more complex systems \cite{LMA+2005}. By contrast, the uniformly second-kind structure of the present formulation suggests improved stability in the time domain. A systematic investigation of the global multi-trace TD-M\"{u}ller formulation is therefore a promising direction for future work.

The rest of this paper is structured as follows. The next section introduces the notation and derives the global multi-trace M\"{u}ller formulation. Section~\ref{sec:discretization} presents a conforming mixed discretization and analyzes its computational cost. A comprehensive set of numerical examples is provided in Section~\ref{sec:results} to validate the accuracy, robustness, and computational efficiency of the proposed method. Finally, concluding remarks are outlined in Section~\ref{sec:conclusion}.

\section{Formulation}
\label{sec:formulation}

\subsection{Problem Setting and Notations}

\begin{figure}
    \centering
    \includegraphics[width=\linewidth]{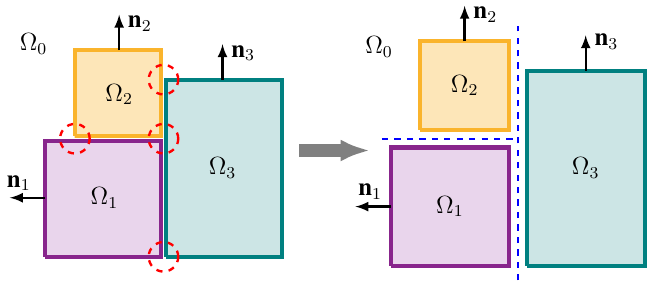}
    \caption{\textit{Left}: A composite object consisting of several components $\Om_k, k = 1, 2, \ldots, N,$ immersed in the background $\Om_0$. This configuration usually gives rise to junctions that are curves where three or more regions meet (dashed red circles). \textit{Right}: Illustration of the global multi-trace approach, in which a conceptual gap filled by the background medium is inserted between adjacent components (dashed blue lines).}
    \label{fig:domain}
\end{figure}

We consider an electromagnetic scatterer composed of $N$  bounded Lipschitz components $\Om_k$, with boundary $\Gm_k := \pa\Om_k$ and outward unit normal vector $\nv_k, k = 1, 2, \ldots, N$. The background region is denoted by $\Om_0 := \R^3 \setminus \bigcup_{k=1}^N \ovl{\Om_k}$ (see Fig.~\ref{fig:domain}, \textit{left}). Each region $\Om_i, i = 0, 1, \ldots, N,$ is filled by a homogeneous, isotropic\footnote{For homogeneous anisotropic media, see \cite{HS2023} for an extension of the M\"{u}ller formulation.}, dielectric material characterized by permittivity $\epsilon_i$ and permeability $\mu_i$. We denote the coefficient matrix associated with the region $\Om_i$
\[
    \beta_i = \diag(\epsilon_i, \mu_i).
\]

For clarity, throughout the paper, the index $i = 0, 1, \ldots, N$ is used for $N + 1$ partitions of $\R^3$ (called regions), whereas the indices $j, k = 1, 2, \ldots, N$ refer to the $N$ components of the scatterer (also called subdomains).

The scatterer is illuminated by time-harmonic incident electric and magnetic fields $(\eb^{inc}, \hb^{inc})$ with angular frequency $\om$. The representation theorem states that any electromagnetic fields $(\eb, \hb)$ satisfying Maxwell's equations in a subdomain $\Om_k, k = 1, 2, \ldots, N,$ can be represented by equivalent surface currents on its boundary $\Gm_k$. More precisely, the Stratton--Chu representation formula associated with the subdomain $\Om_k$ reads \cite{Colton2013}
\begin{align}
    \label{eq:Stratton_Chu}
    & \begin{pmatrix}
        \KK^{(k)}_k & - \eta_k \TT^{(k)}_k \\
        \eta_k^{-1} \TT^{(k)}_k & \KK^{(k)}_k
    \end{pmatrix}
    \begin{pmatrix}
        \eb \times \nv_k \\
        \nv_k \times \hb
    \end{pmatrix} \\
    & \qqqqq =
    \begin{cases}
        (\eb, -\hb)^\transpose \q & \text{if} \q \xb \in \Om_k \\
        (\zrb, \zrb)^\transpose & \text{if} \q \xb \in \R^3 \setminus \ovl{\Om_k},
    \end{cases} \nonumber
\end{align}
where the single- and double-layer potentials are defined by
\begin{align*}
    \TT^{(i)}_k(\mb)(\xb) & = - \iota \kappa_{i} \int_{\Gm_k} \dfrac{e^{-\iota \kappa_i \abs{\xb - \yb}}}{4\pi \abs{\xb - \yb}} \mb(\yb) \ds_{\yb} \\
    & \q + \dfrac{1}{\iota \kappa_i} \gradt \int_{\Gm_k} \dfrac{e^{-\iota \kappa_i \abs{\xb - \yb}}}{4\pi \abs{\xb - \yb}} \divt_\Gm \mb(\yb) \ds_{\yb}, \\
    \KK^{(i)}_{k} (\mb)(\xb) & = \curlt \int_{\Gm_k} \dfrac{e^{-\iota \kappa_i \abs{\xb - \yb}}}{4\pi \abs{\xb - \yb}} \mb(\yb) \ds_{\yb},
\end{align*}
with $i = 0, 1, \ldots, N,$ and $\yb \in \Gm_k, \xb \notin \Gm_k, k = 1, 2, \ldots, N$. Here, $\iota$ denotes the imaginary unit, $\kappa_i = \omega\sqrt{\mu_i \epsilon_i}$ is the wavenumber, and $\eta_i = \sqrt{\mu_i/\epsilon_i}$ is  the impedance coefficient. In the background region $\Omega_0$, the representation formula \eqref{eq:Stratton_Chu} holds for the scattered electromagnetic fields, which satisfy Maxwell's equations together with the radiation condition. The formula \eqref{eq:Stratton_Chu} is also known as the extinction property or the null-field property, meaning that the potential operators associated with one domain reconstruct its fields in the interior but vanish when evaluated outside that domain. The Stratton--Chu representation plays a central role in the construction of BIE formulations for electromagnetic scattering problems. 

Taking the tangential traces of the potential operators yields the single-layer and double-layer BIOs
\[
    T^{(i)}_{jk} := \nv_j \times \TT^{(i)}_k, \qqq
    K^{(i)}_{jk} := \nv_j \times p.v. \, \KK^{(i)}_k,
\]
where $p.v.$ stands for the Cauchy principal value, with $i = 0, 1, \ldots, N,$ and $j, k = 1, 2, \ldots, N$. The following matrix BIO, which has a similar structure to the matrix potential operator in \eqref{eq:Stratton_Chu}, is essential for assembling block matrices in the next sections:
\[
    A^{(i)}_{jk} := 
    \begin{pmatrix}
        K^{(i)}_{jk} & -\eta_i T^{(i)}_{jk} \\
        \eta_i^{-1} T^{(i)}_{jk} & K^{(i)}_{jk}
    \end{pmatrix}.
\]
The operator $A^{(i)}_{jk}$ maps a pair of tangential traces on $\Gm_k$ to a pair of tangential traces on $\Gm_j$ and comprises BIOs defined with respect to the material parameters $(\epsilon_i, \mu_i)$.

For the testing, we also introduce the following local pairing between two pairs $(\pb, \qb)^\transpose$ and $(\jb, \mb)^\transpose$ of tangential traces supported on the boundary $\Gm_k$:
\begin{equation}
    \label{eq:pairing}
    \inprod{
    \begin{pmatrix}
        \pb \\
        \qb
    \end{pmatrix},
    \begin{pmatrix}
        \jb \\
        \mb
    \end{pmatrix}
    }_{\times, k}
    = \int_{\Gm_k} (\nv_k \times \pb) \cdot \jb + (\nv_k \times \qb) \cdot \mb \ds.
\end{equation}
 
\subsection{Global Multi-Trace Method}

The global multi-trace approach seeks the tangential traces of the electric and magnetic fields on the boundary of each component $\Om_k, k = 1, 2, \ldots, N$. On exterior interfaces that are parts of the boundary of $\Om_0$, each trace appears only once. In contrast, on interior interfaces shared by two adjacent components, the traces are duplicated, one for each side of the interface. This distinction can be interpreted in terms of the introduction of an infinitesimally thin gap, so that all $N$ components are treated as floating within the background medium $\Om_0$ and every subdomain contributes its own set of traces to the global system \cite{CHJ2013,CHJ+2014} (see Fig.~\ref{fig:domain}, \textit{right}).
 
Let us denote $\ub := (\ub_1, \ub_2, \ldots, \ub_N)^\transpose,$ where $\ub_k = (\eb \times \nv_k, \nv_k \times \hb)^\transpose$ is the local Cauchy data representing the equivalent magnetic and electric surface current densities on the boundary $\Gm_k$. We assume that the incident fields originate from the background region $\Om_0$ and satisfy Maxwell's equations there. The exterior boundary representation formula is obtained by taking the tangential traces of the Stratton--Chu formula associated with $\Om_0$, from the interior of $\Om_0$ to the boundary of each subdomain $\Om_k, k = 1, 2, \ldots, N$. Since the normal vector on $\Gm_k$ is oriented outward from $\Om_k$, we get
\begin{equation}
    \label{eq:exterior}
    C^{(ext)} \ub = -\ub^{inc},
\end{equation}
where the right-hand vector $\ub^{inc} = (\ub^{inc}_1, \ub^{inc}_2, \ldots, \ub^{inc}_N)^\transpose$, with $\ub^{inc}_k = (\eb^{inc} \times \nv_k, \nv_k \times \hb^{inc})^\transpose$. The $N \times N$ block-matrix operator $C^{(ext)},$ which describes the interaction between subdomains through the background medium, is defined block-wise by \cite{LBG+2022}
\[
    \left[C^{(ext)}\right]_{jk} = A^{(0)}_{jk} - \dfrac{1}{2} I, \qq j, k = 1, 2, \ldots, N,
\]
where $I$ is a so-called geometric identity, as introduced in \cite{LBG+2022}. This means, in theory, the identity operator is included only when $\Gm_j$ and $\Gm_k$ share a common interface. In practical boundary-element implementations, however, it is unnecessary to distinguish explicitly between overlapping and non-overlapping interfaces when deriving the formulation. The geometric support of the basis functions implicitly enforces this distinction: when the two boundaries are disjoint, the corresponding discretized operator simply reduces to a zero matrix.

The interior boundary representation formula is derived analogously by taking the tangential traces of the Stratton--Chu formula \eqref{eq:Stratton_Chu} associated with each $\Om_k, k = 1, 2, \ldots, N$, from the interior of $\Om_k$ to its boundary $\Gm_k$, yielding
\begin{equation}
    \label{eq:interior}
    C^{(int)} \ub = \zrb,
\end{equation}
where the block-diagonal operator $C^{(int)}$ is defined by
\[
    \left[C^{(int)}\right]_{kk} = A^{(k)}_{kk} + \dfrac{1}{2} I, \qq k = 1, 2, \ldots, N.
\]

The exterior and interior representation formulas \eqref{eq:exterior} and \eqref{eq:interior} alone do not provide a unique solution for the Cauchy data $\ub$. They must be combined into a single formulation that enforces the transmission conditions appropriately. The most widely used approach is the global multi-trace (MT-) PMCHWT equation, obtained by summing the two representations \cite{CH2012}
\begin{equation}
    \label{eq:PMCHWT}
    P \ub := \paren{C^{(ext)} + C^{(int)}} \ub = -\ub^{inc}.
\end{equation}
The blocks of the operator $P$ are given by
\[
    \left[ P \right]_{kk} = A^{(0)}_{kk} + A^{(k)}_{kk}, \qq \left[ P \right]_{jk} = A^{(0)}_{jk} - \dfrac{1}{2} I \q (j \neq k).
\]
The MT-PMCHWT formulation \eqref{eq:PMCHWT} is a first-kind BIE. In the next section, we derive a second-kind global multi-trace formulation based on the M\"{u}ller equation.

\subsection{Multi-Trace M\"{u}ller Formulation}

In contrast to the PMCHWT equation, the M\"{u}ller formulation is constructed as a weighted linear combination of the exterior and interior representation formulas, with the weights chosen in terms of the material coefficients, such that the hyper-singularities cancel \cite{Muller1969}. Following this principle, we consider the formulation
\begin{equation}
    \label{eq:Mueller_w}
    Q \ub := \paren{B C^{(int)} - B_0 C^{(ext)}} \ub = B_0 \ub^{inc},
\end{equation}
where the global coefficient matrices
\[
    B_0 = \diag(\beta_0, \beta_0, \ldots, \beta_0), \qq B = \diag(\beta_1, \beta_2, \ldots, \beta_N).
\] 
The operator $Q$ can be determined block-wise by
\begin{align*}
    \left[Q\right]_{kk} & = \dfrac{1}{2} \paren{\beta_k + \beta_0} I + \beta_k A^{(k)}_{kk} - \beta_0 A^{(0)}_{kk}, \\
    \left[Q\right]_{jk} & = \dfrac{1}{2} \beta_0 I - \beta_0 A^{(0)}_{jk} \qqqq (j \neq k).
\end{align*}
In the diagonal block $\left[Q\right]_{kk}$, the hyper-singular contributions come from the single-layer operators $T_{kk}^{(k)}$ and $T_{kk}^{(0)}$ in $A^{(k)}_{kk}$ and $A^{(0)}_{kk}$, associated with and weighed by the material parameters $\beta_k$ and $\beta_0$ of the regions $\Om_k$ and $\Om_0$, respectively. These contributions have the same prefactor, namely 1. Their difference has the kernel
\begin{align*}
    \dfrac{e^{-\iota \kappa_k \abs{\xb - \yb}} - e^{-\iota \kappa_0 \abs{\xb - \yb}}}{4 \, \iota \pi \om \abs{\xb - \yb}} & \sim \dfrac{\sqrt{\mu_0\eps_0} - \sqrt{\mu_k \eps_k}}{4\pi} \\
    & \,\,\, + \dfrac{\mu_k \eps_k - \mu_0\eps_0}{8\pi} \iota \omega \abs{\xb - \yb},
\end{align*}
as $\omega \abs{\xb - \yb} \to 0$, and is thus regular with respect to both the frequency and spatial variables. Moreover, the leading term on the right-hand side is constant and vanishes upon application of the exterior gradient, further increasing the regularity of the kernel. Consequently, the diagonal blocks of $Q$ contain no hyper-singularities, and are therefore second-kind operators. In contrast, the off-diagonal blocks exhibit hyper-singularity, which can result in ill-conditioned linear systems when discretized on dense meshes and at low frequencies. This inconsistency has contributed to the long-standing perception that the singularity-cancellation mechanism underlying the M\"{u}ller equation breaks down in composite configurations \cite{TY2006}.

The authors of \cite{CHS2017} resolved this issue by leveraging the so-called multi-potential representation formula, which exploits the inherent continuity of the single-trace solutions to regularize the system and obtain a fully consistent second-kind equation. That technique, however, is not applicable in our global multi-trace framework, as the background region is treated differently from the $N$ subdomains of the scatterer. In particular, on exterior interfaces, we consider only one pair of trace unknowns instead of two, and therefore the multi-potential formula does not hold on those interfaces.

In this contribution, we rely on the Stratton--Chu representation formula to regularize the M\"{u}ller formulation. More specifically, for each pair $j \neq k$, we augment the off-diagonal block $(j, k)$ of the interior representation formula with the tangential traces of the Stratton--Chu potentials associated with the subdomain $\Om_k$, evaluated from the exterior of $\Om_k$ onto the boundary $\Gm_j$. The resulting BIO is given by 
\begin{equation}
    \label{eq:augment_jk}
    A^{(k)}_{jk} - \dfrac{1}{2} I.
\end{equation}
To ensure the cancellation of all hyper-singularities in the resulting global multi-trace (MT-) M\"{u}ller formulation, the coefficients of the exterior and interior representation operators must also be consistent. Specifically, all hyper-singular terms in the weighted exterior operator $B_0 C^{(ext)}$ appear with prefactor $1$, whereas those in the augmenting operator $A^{(k)}_{jk}$, after weighing by $\beta_j$ in $B$, carry prefactor $\beta_j \beta_k^{-1}$. Therefore, the augmenting $(j, k)$ block \eqref{eq:augment_jk} must be rescaled by the material ratio $\beta_j^{-1} \beta_k$ so that its hyper-singular contribution also appears with prefactor 1 and will exactly cancel the leading singular part of the exterior contribution. The modified interior representation operator is then defined as
\begin{align}
    \left[\wh{C}^{(int)}\right]_{kk} & = \left[C^{(int)}\right]_{kk} = A^{(k)}_{kk} + \dfrac{1}{2} I, \notag \\
    \left[\wh{C}^{(int)}\right]_{jk} & = \beta_j^{-1} \beta_k \paren{A^{(k)}_{jk} - \dfrac{1}{2} I} \q (j \neq k). \label{eq:augmented_op}
\end{align}
Due to the extinction property, the augmenting operators vanish when acting on the corresponding field traces; that is, for all $j, k = 1, 2, \ldots, N$ and $j \neq k$
\[
    \left[\wh{C}^{(int)}\right]_{jk} \ub_k = \zrb.
\]
Thus, the modified interior representation formula reads
\begin{equation}
    \label{eq:mod_interior}
    \wh{C}^{(int)} \ub = \zrb.
\end{equation}
Combining \eqref{eq:mod_interior} with the exterior representation formula \eqref{eq:exterior} yields the MT-M\"{u}ller formulation
\begin{equation}
    \label{eq:Mueller}
    M \ub := \paren{B \wh{C}^{(int)} - B_0 C^{(ext)}} \ub = B_0 \ub^{inc},
\end{equation}
where the block entries of $M$ are given by
\begin{align*}
    \left[M\right]_{kk} & = \dfrac{1}{2} \paren{\beta_k + \beta_0} I + \beta_k A^{(k)}_{kk} - \beta_0 A^{(0)}_{kk}, \\
    \left[M\right]_{jk} & = \dfrac{1}{2} \paren{\beta_0 - \beta_k} I + \beta_k A^{(k)}_{jk} - \beta_0 A^{(0)}_{jk} \q\, (j \neq k).
\end{align*}
It is clear that every block of the operator $M$ exhibits second-kind structure.

In the case of a single homogeneous dielectric (i.e., $N=1$), the formulation \eqref{eq:Mueller} reduces to the classical M\"{u}ller equation. If a component $\Om_k$ is filled with the background medium, then all off-diagonal entries in the $k$-th column of the operator $M$ vanish, while the diagonal block becomes $\beta_0 I$. When all components coincide with the background medium (i.e., $B = B_0$), the formulation \eqref{eq:Mueller} degenerates to the trivial identity
\begin{equation}
    \label{eq:trivial}
    B_0 \ub = B_0 \ub^{inc}.
\end{equation}

\section{Discretization}
\label{sec:discretization}

Let the boundary $\Gm_k$ of subdomain $\Om_k, k = 1, 2, \ldots, N,$ be discretized by a triangulation $\Gm_{h, k}$ with average mesh size $h$, consisting of $N_{e, k}$ edges. The global mesh is defined as the disjoint union 
\[
    \Gm_h := \bigsqcup_{k = 1}^N \Gm_{h, k},
\]
which contains a total of $N_e := \sum_{k = 1}^N N_{e, k}$ edges. Notably, $\Gamma_{h, k}$ and $\Gamma_{h, j}$ are not required to be mutually conformal on the shared interface $\Gamma_j \cap \Gamma_k$. This is a key advantage of the multi-trace framework, as it allows the subdomain boundaries to be discretized independently, thereby enabling the efficient numerical treatment of geometrically complex, large-scale, and multi-scale electromagnetic problems. This flexibility is particularly valuable in practical design and simulation settings, where enforcing global mesh conformity may be unnecessarily restrictive or computationally expensive. Finally, although the notation $\Gm_h$ resembles the notation used for the boundaries $\Gm_k$, the subscript $h$ denotes discrete quantities and should not be interpreted as an additional index.

\subsection{Mixed Discretization}
 
The MT-M\"{u}ller equation \eqref{eq:Mueller} is discretized using a conforming Petrov--Galerkin (mixed) scheme, employing Rao--Wilton--Glisson (RWG) basis functions \cite{RWG1982} for the trial space and BC functions for the test space. The multi-trace trial space is defined as the direct product
\[
    \mathbb{RT}(\Gm_h) := \prod_{k=1}^N \mathbf{RT}(\Gm_{h, k}),  
\]
where $\mathbf{RT}(\Gm_{h, k}) := \text{RT}(\Gm_{h, k}) \times \text{RT}(\Gm_{h, k})$, and $\text{RT}(\Gm_{h, k})$ is the boundary element space spanned by the $N_{e, k}$ RWG functions defined on the mesh $\Gm_{h, k}$ of $\Gm_k$. Analogously, the test space is defined as
\[
    \mathbb{BC}(\Gm_h) := \prod_{k=1}^N \mathbf{BC}(\Gm_{h, k}),  
\]
where $\mathbf{BC}(\Gm_{h, k}) := \text{BC}(\Gm_{h, k}) \times \text{BC}(\Gm_{h, k})$, with $\text{BC}(\Gm_{h, k})$ spanned by the BC basis functions on $\Gm_{h, k}$. This space is a discrete subspace of the multi-trace energy space and is locally dual to the trial space, meaning that the components associated with each subdomain are mutually dual.

The trace unknown $\ub$ is approximated by its expansion in $\mathbb{RT}(\Gm_h)$ and the corresponding coefficient vector is denoted by $\uv \in \C^{2N_e}$. The mixed discrete formulation of \eqref{eq:Mueller} is then given by
\begin{equation}
    \label{eq:discretization}
    \M \, \uv = \uv^{inc},
\end{equation}
where $\M$ is the matrix representation of the bilinear form $\inprod{\cdot, M \cdot}_{\times}$ on $\mathbb{BC}(\Gm_h) \times \mathbb{RT}(\Gm_h)$ with respect to its standard basis, and $\uv^{inc}$ is the coefficient vector representing the linear functional $\inprod{\cdot, B_0 \ub^{inc}}_{\times}$ on $\mathbb{BC}(\Gm_h)$. Here, the global pairing $\inprod{\cdot, \cdot}_{\times}$ is defined by
\[
    \inprod{\fv, \gv}_\times := \sum_{k=1}^N \inprod{\fv_k, \gv_k}_{\times, k}.
\]

\subsection{Computational Cost Analysis}

The mixed discretization \eqref{eq:discretization} yields linear systems that remain well conditioned across a wide range of frequencies and mesh densities. Consequently, iterative Krylov solvers converge rapidly without requiring any additional stabilization, thereby avoiding the cost associated with extra matrix-vector evaluations. Numerical results consistently demonstrate this advantage: the MT-M\"{u}ller formulation achieves significantly shorter solving times than both the MT-PMCHWT equation and its CP variants. Although CP MT-PMCHWT formulations often exhibit a mesh-independent iteration count, their per-iteration matrix-vector multiplication costs grow quickly with mesh refinement and scale up the total solving time. This advantage of the MT-M\"{u}ller formulation becomes even more pronounced when multiple right-hand sides are involved.

Despite these benefits, the MT-M\"{u}ller equation presents two noteworthy drawbacks. The first stems from the mixed discretization of second-kind BIEs. Constructing BC test functions necessitates a barycentric refinement of the surface mesh, which increases the assembly time of the MT-M\"{u}ller operator by roughly a factor of six compared with testing on the original mesh. It is worth noting, however, that efficient Calder\'{o}n preconditioners for the PMCHWT formulation also rely on BC or CW functions and incur an even larger assembly cost, about 36 times that of the non-refined discretization. Refinement-free Calder\'{o}n strategies \cite{AAE2019} or OSRC preconditioners \cite{FB2023} can mitigate this burden, but they require non-trivial matrix manipulations or exhibit inconsistent performance across different frequency ranges.

The second drawback concerns the increased number of BIOs in the MT-M\"{u}ller formulation. In the MT-PMCHWT equation \eqref{eq:PMCHWT}, each diagonal block contains two single-layer and two double-layer operators, while each off-diagonal block consists of only one identity operator, one single-layer operator, and one double-layer operator (counting only unique operators due to symmetry). In contrast, every block of the MT-M\"{u}ller operator in \eqref{eq:Mueller} includes one identity, two single-layer, and two double-layer operators. Although identity operators are local and can be assembled efficiently into sparse matrices, the additional non-local operators substantially increase the assembly time, especially for configurations with many scattering components.

In summary, the MT-M\"{u}ller formulation offers superior solving times but higher assembly costs. For moderate numbers of sources and scattering components, its overall computational cost is broadly comparable to that of the CP MT-PMCHWT equation. When the number of excitations increases, the MT-M\"{u}ller equation becomes more efficient. Conversely, for large numbers of subdomains, the CP MT-PMCHWT approach is more favorable. The latter situation is typically not a limitation in practice, as hybrid finite-element--boundary-element coupling strategies generally outperform both formulations in such scenarios and thus constitute the preferred approach. 

Finally, for a fixed number of components, the MT-M\"{u}ller and MT-PMCHWT formulations exhibit similar asymptotic complexities when accelerated with hierarchical compression schemes such as Adaptive Cross Approximation (ACA) \cite{Bebendorf2000,TA2024}, H-matrix method \cite{BGH2003}, or Fast Multipole Methods (FMMs) \cite{Rokhlin1990,CRW1993}.

\section{Numerical Results}
\label{sec:results}

\begin{figure*}[!t]
    \centering
    \includegraphics[trim={1.2cm 1.2cm 1.2cm 1.2cm}, clip, width=0.31\linewidth]{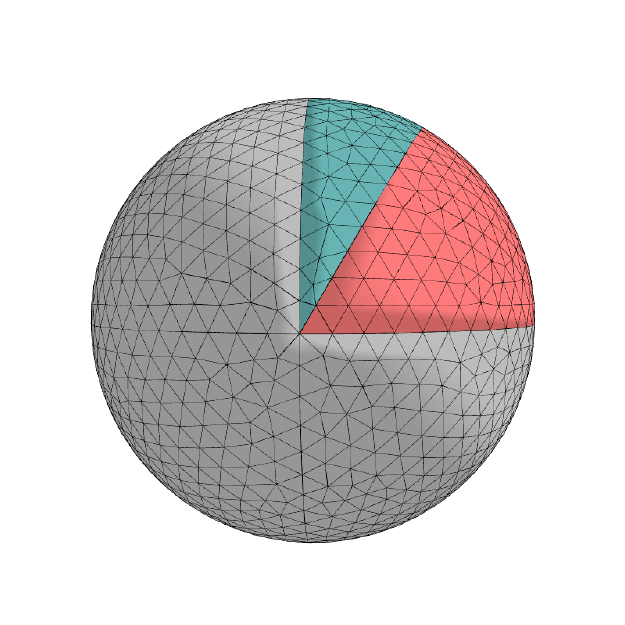}
    \hfill
    \includegraphics[trim={1.4cm 2.3cm 2.2cm 3cm}, clip, width=0.35\linewidth]{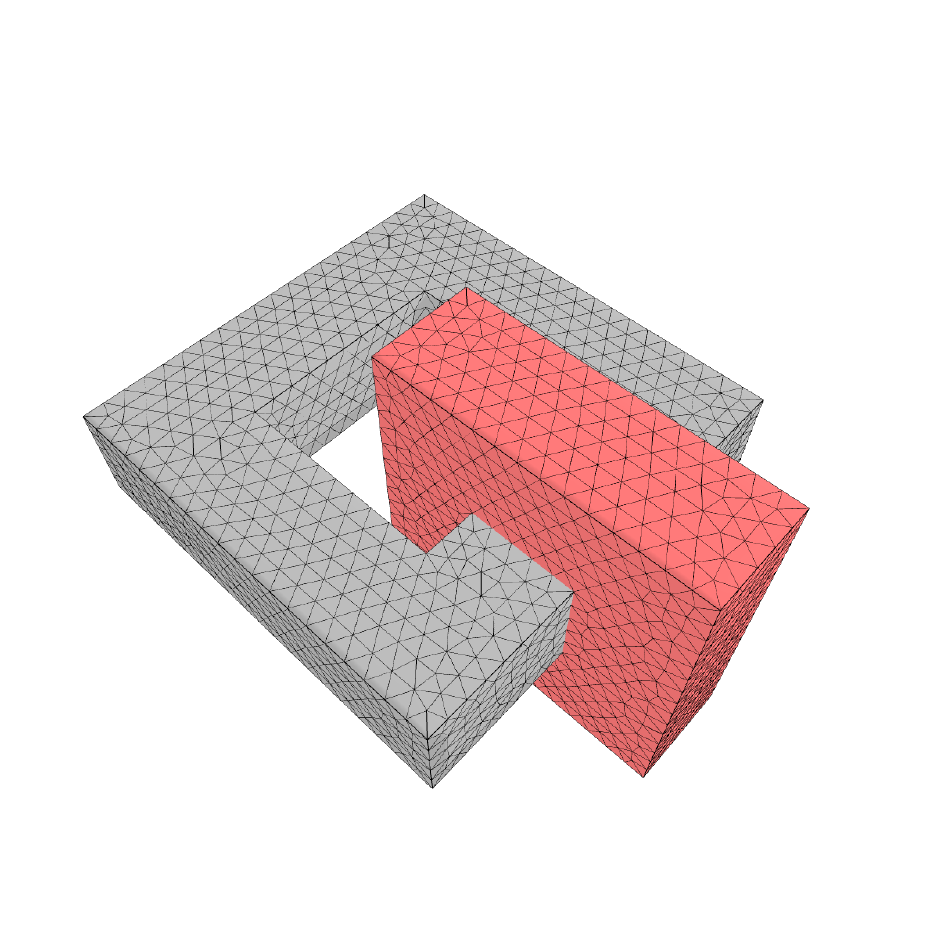}
    \hfill
    \includegraphics[trim={1.4cm 1.4cm 1.4cm 1.5cm}, clip, width=0.29\linewidth]{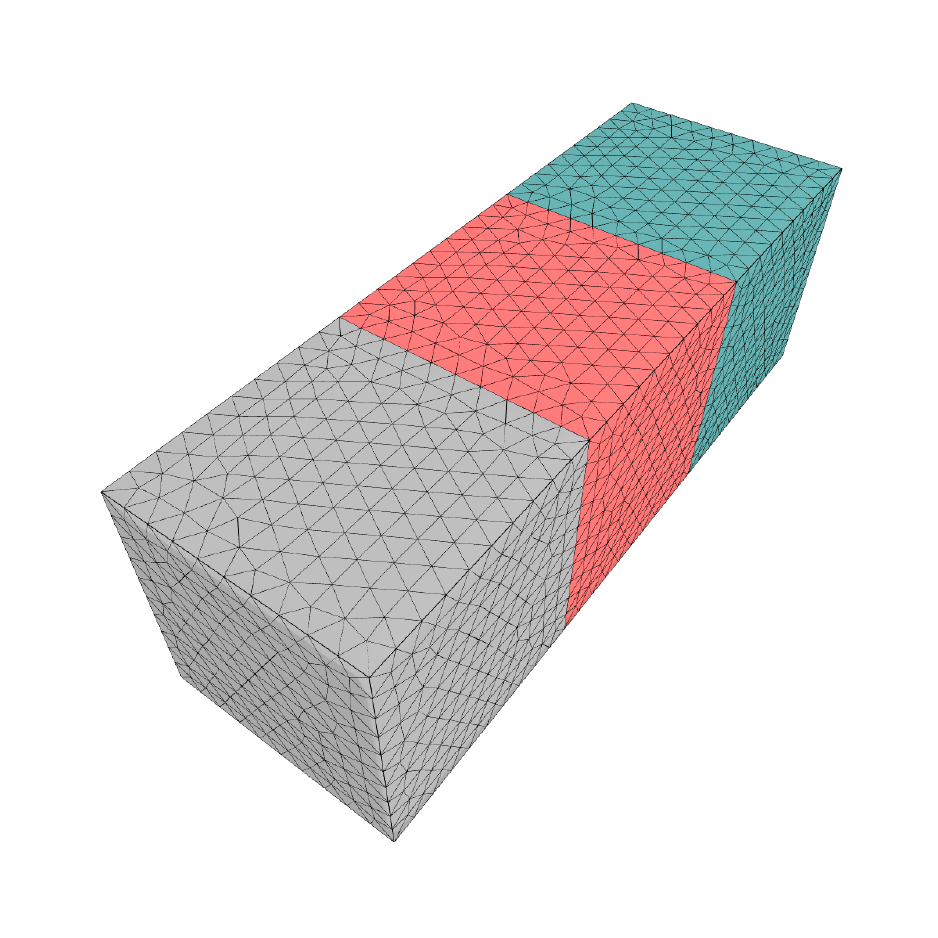}
    \caption{Geometries used in numerical experiments. From left to right: a three-component unit sphere in which one component occupies three quadrants, the second occupies two thirds of the remaining quadrant, and the third occupies the rest; a configuration of two fused square-section tori featuring a fully occluded cavity; and a domain constructed by stacking multiple unit cubes (illustrated here using three cubes).}
    \label{fig:geometries}
\end{figure*}

In this section, we present numerical experiments that validate the accuracy, stability, and computational efficiency of the proposed MT-M\"{u}ller formulation. The following three geometries with distinct configurations and material contrasts are considered (see Fig.~\ref{fig:geometries}):
\begin{itemize}
    \item[(a)] A sphere of radius $1\mathrm{m}$, partitioned into three components: one occupying three quadrants, a second occupying two thirds of the remaining quadrant, and a third occupying the rest. Apart from the high-contrast experiments, all three components are filled with the same material characterized by $(\epsilon, \mu) = (3\epsilon_0, \mu_0)$.
    \item[(b)] A dual-torus geometry formed by the intersection of two square-section tori of dimensions $1\mathrm{m} \times 1\mathrm{m} \times 0.25\mathrm{m}$ and $0.75\mathrm{m} \times 0.25 \mathrm{m} \times 0.75\mathrm{m}$, such that one cavity is completely occluded. The two subdomains are filled with materials $(\epsilon_1, \mu_1) = (2\epsilon_0, \mu_0)$ and $(\epsilon_2, \mu_2) = (4\epsilon_0, \mu_0)$.
    \item[(c)] A vertical stack (along the $z$-axis) of $N$ cubes of side length $1\mathrm{m}$. The $k$-th cube (counted from bottom to top) is filled with material $(\epsilon_k, \mu_k) = \paren{(k+1)\epsilon_0, \mu_0}$.
\end{itemize}
It is worth noting that all these geometries contain junctions where either three subdomains meet or two subdomains meet the background region. Throughout this section, the background medium is free space and, unless otherwise specified, the incident fields are the plane wave of the form
\begin{align*}
    \eb^{inc}(\xb) & = \hat{\xb} \exp(-\iota \kappa_0 \hat{\zb} \cdot \xb), \\
    \hb^{inc}(\xb) & = -\dfrac{1}{\iota \kappa_0 \eta_0} \curlt \, \eb^{inc}(\xb),
\end{align*}
where $\hat{\xb}$ and $\hat{\zb}$ denote the unit vectors along the $x$-axis and $z$-axis, respectively.

\begin{figure*}[!t]
    \centering
    \includegraphics[trim={1.5cm 1.4cm 1cm 1.5cm}, clip, width=0.28\linewidth]{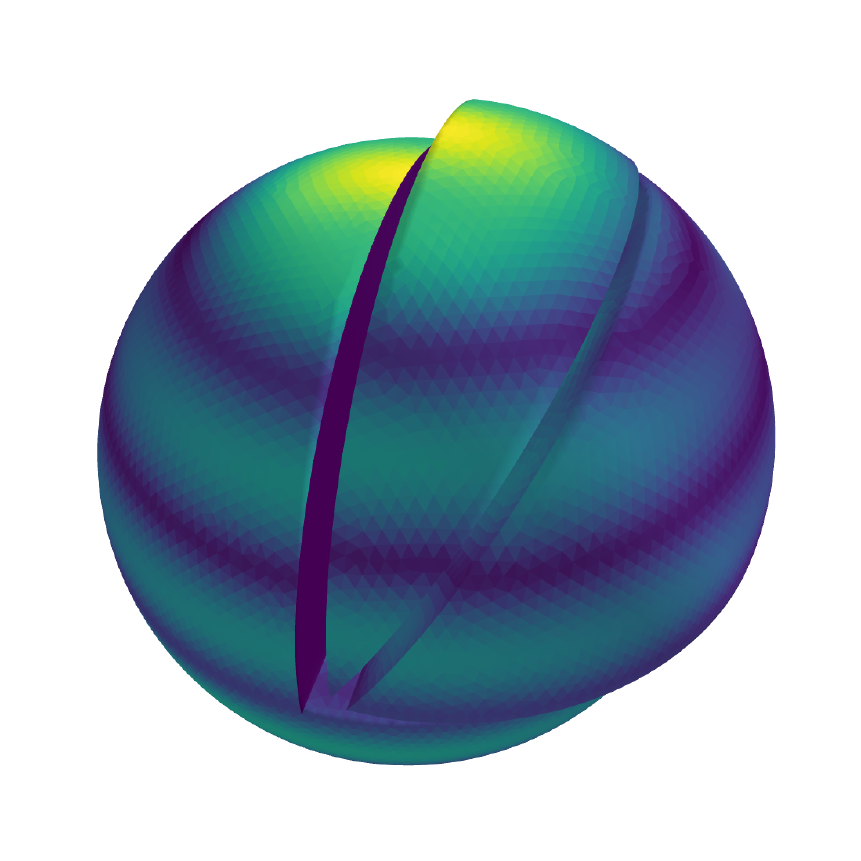}
    \hfill 
    \includegraphics[trim={1.2cm 1.4cm 1cm 1.5cm}, clip, width=0.37\linewidth]{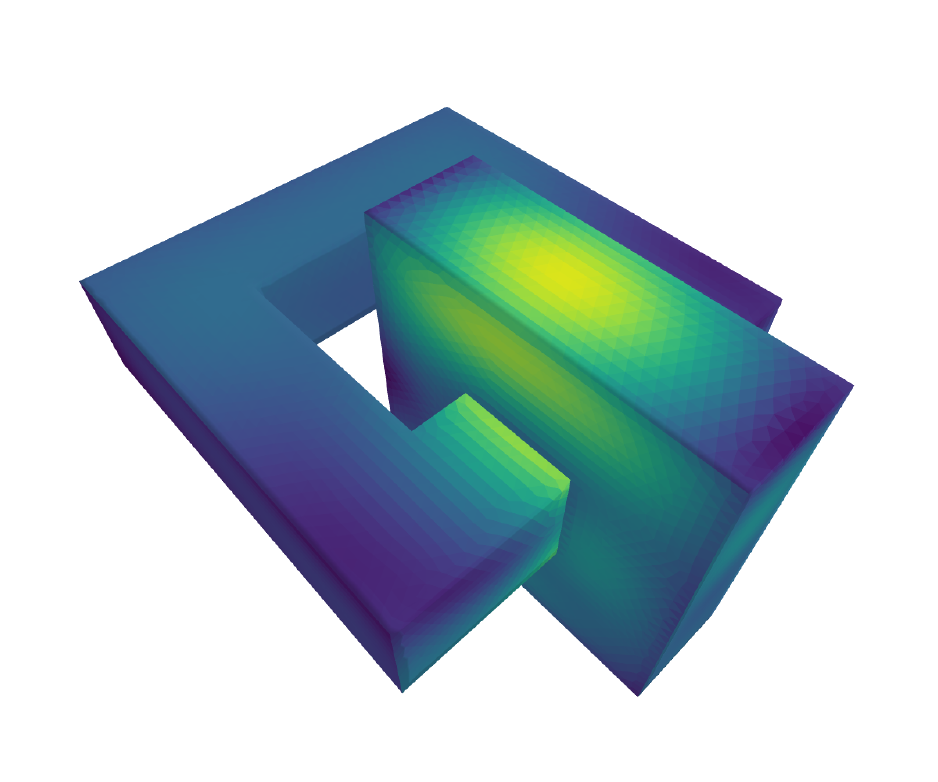}
    \hfill
    \includegraphics[trim={1.6cm 1.2cm 0.4cm 1.5cm}, clip, width=0.32\linewidth]{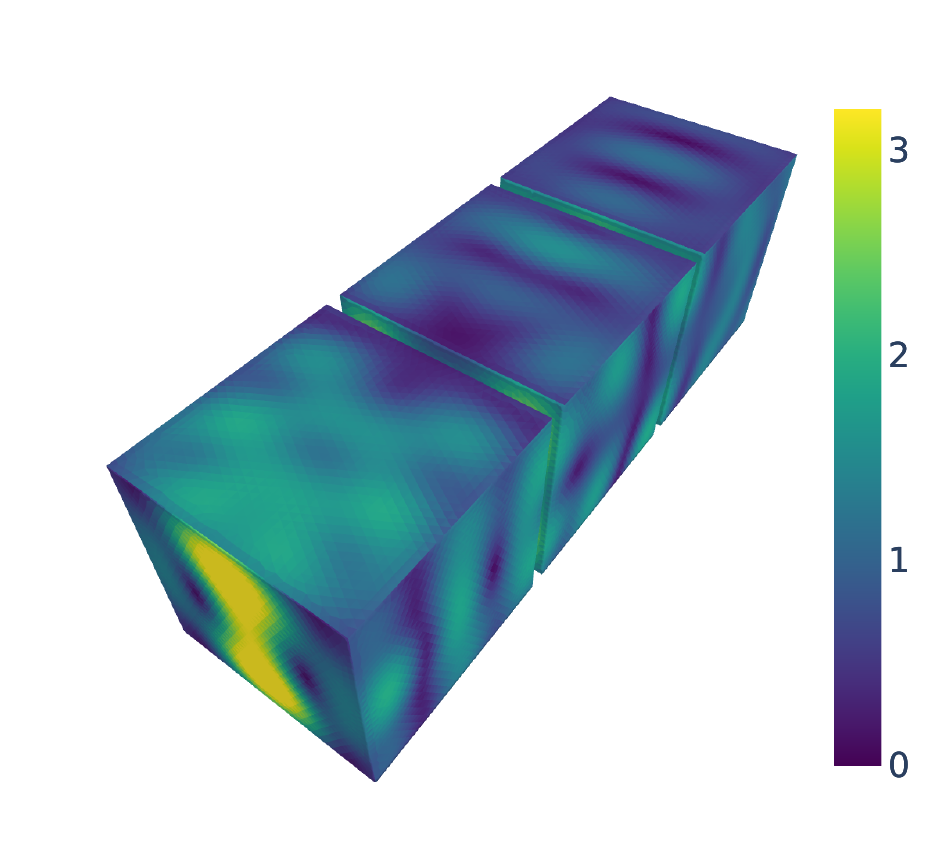}
    \caption{Magnitude of the electric surface current density on the boundaries of the three geometries, excited by a plane wave with wavenumber $\kappa_0 = 6\mathrm{m^{-1}}$, computed using the global multi-trace M\"{u}ller equation. For visualization purposes, the subdomains of the sphere (leftmost) and the cubes (rightmost) are slightly separated. The solutions exhibit continuity across the interfaces between adjacent materials.}
    \label{fig:facecurrent}
\end{figure*}

Fig.~\ref{fig:facecurrent} displays the magnitude of the electric current density $\nv \times \hb$ on the boundaries of the three domains, computed using the mixed MT-M\"{u}ller formulation \eqref{eq:discretization}. In all cases, the structures are excited by a plane wave with wavenumber $\kappa_0 = 6\mathrm{m^{-1}}$. The computed currents exhibit the expected continuity across interfaces between adjacent materials.

\subsection{Extinction Property}

\begin{figure}[http]
    \centering
    \includegraphics[width=\linewidth]{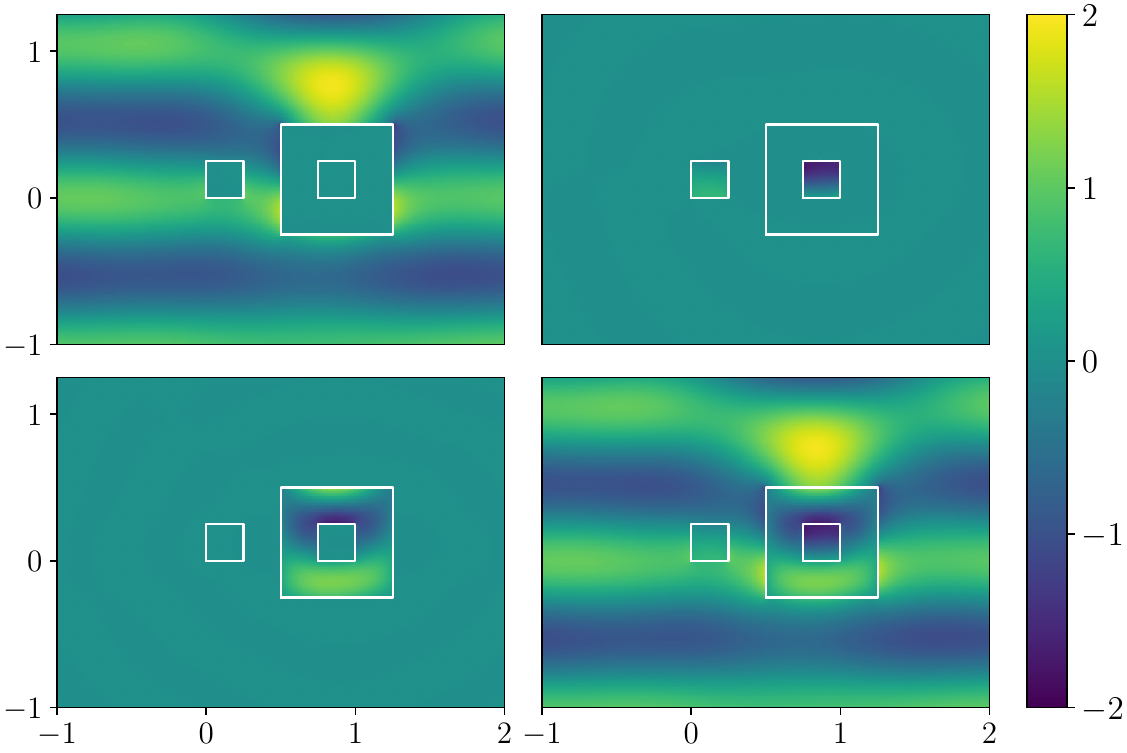}
    \caption{Real part of the $y$-component of the electric field reconstructed from the Cauchy data on the interfaces of the dual-torus domain using the Stratton--Chu representation potentials, evaluated on a grid at $y = 0.5 \mathrm{m}$. From top to bottom and left to right, the subplots show the fields reconstructed using the potentials associated with the background region, with the two subdomains, and finally the total electric field. The surface currents fulfill the extinction property and therefore constitute a Maxwell solution. The total field also exhibits the continuity across the material interfaces.}
    \label{fig:nearfield}
\end{figure}

In this section, we validate the correctness of the proposed MT-M\"{u}ller formulation by examining the extinction property \eqref{eq:Stratton_Chu} and the continuity of the derived near fields for the simulation shown in Fig.~\ref{fig:facecurrent} (middle). The electric and magnetic fields are reconstructed both inside and outside the dual-torus structure using the Stratton--Chu representation potentials. Fig.~\ref{fig:nearfield} presents the reconstructed electric field evaluated on a grid at $y = 0.5 \mathrm{m}$. The first three subplots correspond to the fields generated using the potentials associated with the material coefficients of the background region and the two interior subdomains, while the final subplot shows the total electric field obtained by summing the contributions from all regions.

The reconstructed fields confirm that the computed surface currents satisfy the extinction property, therefore constitute a solution to Maxwell's equations. Moreover, the total field exhibits the continuity across material interfaces. These observations demonstrate that the proposed MT-M\"{u}ller formulation correctly enforces transmission conditions within the discretization error and provides an accurate boundary representation of the scattering problem.

\subsection{Far-Field Evaluation}

\begin{figure*}[!t]
    \centering
    \includegraphics[trim={0.4cm 0cm 0.4cm 0cm}, clip, width=\linewidth]{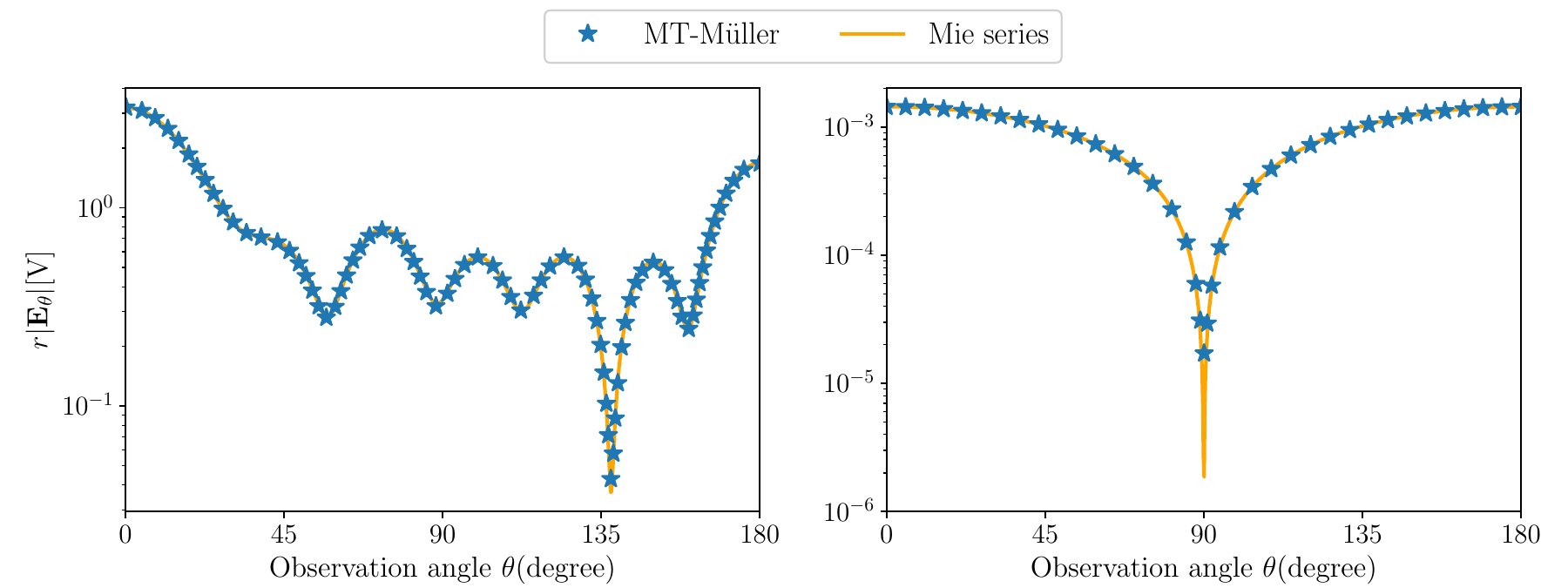}
    \caption{Electric far field scattered by the homogeneous unit sphere for incident plane waves with wavenumbers $\kappa_0 = 6 \mathrm{m^{-1}}$ \textit{(left)} and $\kappa_0 = 0.06\mathrm{m^{-1}}$ \textit{(right)}. The far fields computed using the MT-M\"{u}ller formulation show excellent agreement with those obtained from the Mie series.}
    \label{fig:farfield}
\end{figure*}

To assess the accuracy of secondary quantities produced by the proposed method, we first evaluate the far-field patterns obtained from the surface current densities computed on the homogeneous sphere and compare them with the analytical Mie-series solutions \cite{HRA2023}. Fig.~\ref{fig:farfield} illustrates the scattered electric far fields for two incident plane waves with background wavenumbers $\kappa_0 = 6\mathrm{m^{-1}}$ and $\kappa_0 = 0.06\mathrm{m^{-1}}$. In both cases, the MT-M\"{u}ller formulation reproduces the exact far fields with high fidelity.

\begin{figure}[!t]
    \centering
    \includegraphics[trim={1.5cm 1.3cm 1.3cm 1.3cm}, clip, width=\linewidth]{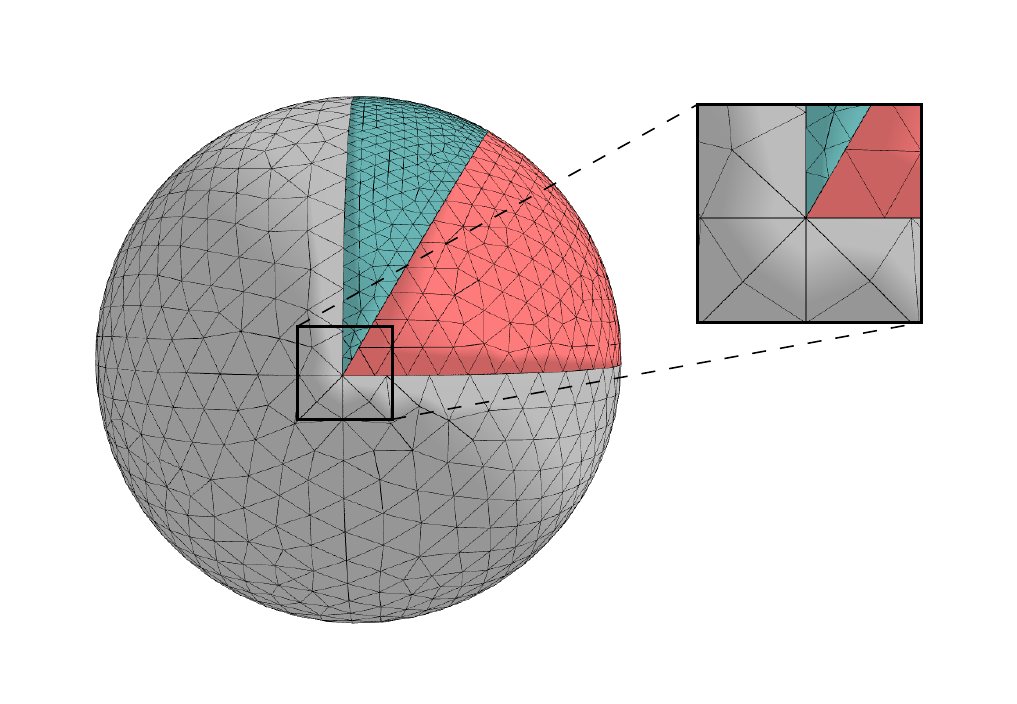}
    \caption{Non-conformal discretization of the unit sphere in the high-contrast heterogeneous configuration. Each subdomain boundary is triangulated independently with mesh size $h = 0.05\lambda$, where $\lambda$ denotes the wavelength in the corresponding subdomain at an incident wavenumber  $\kappa_0 = 1\mathrm{m}^{-1}$.}
    \label{fig:nonconformalmesh}
\end{figure}

Next, we assess the accuracy of the mixed MT-M\"{u}ller formulation in the presence of high-contrast materials. To this end, we consider scattering by the three-component sphere, whose subdomains are filled with materials characterized by $(\eps_1, \eps_2, \eps_3) = (7\eps_0, 14 \eps_0, 42\eps_0)$ and $\mu_1 = \mu_2 = \mu_3 = \mu_0$. In this setting, a non-conformal mesh is a natural choice, as it allows the oscillatory field behavior in each subdomain to be resolved consistently. Fig.~\ref{fig:nonconformalmesh} shows a non-conformal discretization of the sphere, in which the boundary of each subdomain is triangulated with a mesh size $h$ chosen as one twentieth of the wavelength of the corresponding interior field for an incident wavenumber $\kappa_0 = 1\mathrm{m}^{-1}$.

Since non-conformal interfaces prevent the direct application of standard assembly procedures, slight modifications to the BIO quadrature routine are required to evaluate the interaction matrices accurately. In this work, we adopt the simple yet efficient quadrature strategy proposed in \cite{Muenger2024,MC2026}. This strategy performs on-the-fly local refinement for each pair of triangles encountered during the assembly whenever the triangles overlap or share a non-negligible part of their edges. The resulting locally refined meshes are conformal, allowing standard quadrature techniques, such as the Sauter--Schwab singularity cancellation and the Wilton singularity extraction, to be applied. Alternative approaches for electromagnetic scattering on non-conformal meshes were proposed in \cite{SUR2018}, based instead on non-conforming discretizations using the facet-based monopolar-RWG set.

\begin{figure*}[!t]
    \centering
    \includegraphics[trim={0.4cm 0cm 0.4cm 0cm}, clip, width=\linewidth]{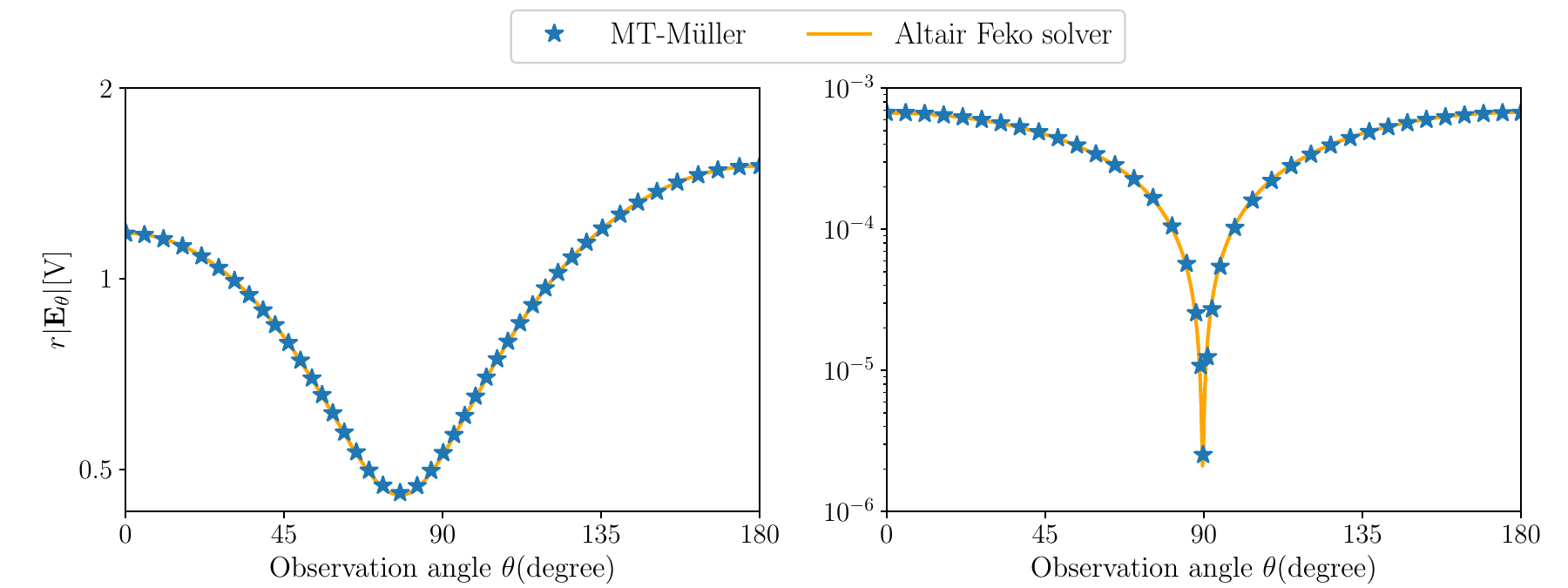}
    \caption{Scattered electric far field of the heterogeneous sphere for incident plane waves with wavenumbers $\kappa_0 = 1\mathrm{m^{-1}}$ \textit{(left)} and $\kappa_0 = 0.03\mathrm{m^{-1}}$ \textit{(right)}. The far fields computed using the MT-M\"{u}ller formulation agree very well with the results obtained from the commercial software Altair Feko.}
    \label{fig:farfield_nonconformal}
\end{figure*}

Fig.~\ref{fig:farfield_nonconformal} presents the scattered electric far fields computed from the MT-M\"{u}ller solutions for two incident plane waves with $\kappa_0 = 1\mathrm{m}^{-1}$ and $\kappa_0 = 0.03\mathrm{m}^{-1}$. Since Mie-series solutions are not available for this heterogeneous configuration, we compare the results with reference solutions obtained using the commercial electromagnetic solver Altair Feko 2026. The two far-field solutions agree very well. These results demonstrate that the mixed discretization of the MT-M\"{u}ller formulation provides accurate far-field predictions for both low- and moderately high-contrast dielectric scatterers in the moderate- and high-frequency regimes.

In the low-frequency limit, the solenoidal and irrotational components of the trace solution may scale differently, leading to numerical cancellation when both components are stored within a single floating-point variable. This cancellation results in inaccurate far-field computations \cite{ADC+2021,LCA+2024}. As shown in \cite{BCA+2014}, the mixed M\"{u}ller formulation produces traces with the correct scaling up to frequencies that solely depend on quadrature and machine precisions, thereby allowing for the introduction of a stabilization scheme. However, a detailed investigation of low-frequency stabilization techniques for the M\"{u}ller and MT-M\"{u}ller formulations lies beyond the scope of the present work and will be addressed in a forthcoming study.

\subsection{Conditioning}

The defining advantage of second-kind BIEs over their first-kind counterparts lies in their favorable conditioning properties. To substantiate this advantage for the proposed MT-M\"{u}ller formulation, we investigate its numerical conditioning by examining both the condition number of the resulting system matrix $\M$ and the number of GMRES iterations required to reach a prescribed convergence tolerance, over a broad range of mesh resolutions and frequencies.

The dual-torus domain is chosen as a particularly revealing test geometry due to its non-trivial topological structure. Each subdomain is itself a torus with genus $g = 1$, while the boundary of the complement region $\Om_0$ is also multiply connected with genus $g = 1$ (as a consequence of the occlusion of an interior cavity), leading to a strongly coupled topology. It is well known that the static interior and exterior magnetic field integral operators $\frac{1}{2} I \pm \widetilde{K}$, with $\widetilde{K}$ denoting the zero-frequency double-layer BIO, possess non-trivial nullspaces on toroidal surfaces. As the frequency tends to zero, these nullspaces manifest numerically as a growth in the condition number of the corresponding discretized operators \cite{CAO+2009a}. The intricate topology of the dual-torus geometry suggests that the proposed MT-M\"{u}ller operator may exhibit a similar growth in its condition number when the frequency approaches zero, since all of its blocks involve either the weighted summation or subtraction of two magnetic field integral operators on their diagonals.

For comparison, we recall the standard discretization of the MT-PMCHWT equation \eqref{eq:PMCHWT}
\begin{equation}
    \label{eq:dis_PMCHWT}
    \Pbb \, \uv = \wh{\uv}^{inc},
\end{equation}
where $\Pbb$ is the matrix representation of the bilinear form $\inprod{\cdot, P \cdot}_\times$ on $\mathbb{RT}(\Gm_h) \times \mathbb{RT}(\Gm_h)$ and $\wh{\uv}^{inc}$ is the vector corresponding to the linear functional $\inprod{\cdot, -\ub^{inc}}_{\times}$ on $\mathbb{RT}(\Gm_h)$. 

\begin{figure*}[!t]
    \centering
    \includegraphics[trim={0.5cm 0cm 0.5cm 0cm}, clip, width=\linewidth]{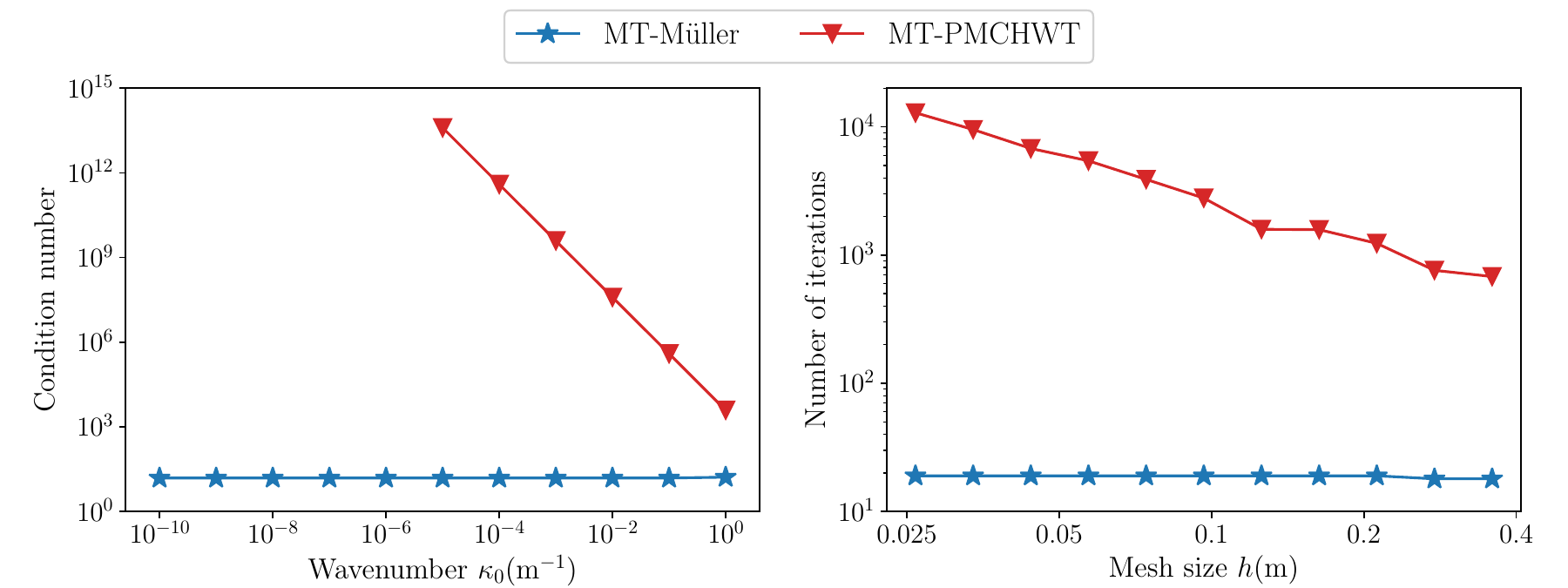}
    \caption{Conditioning behavior of the MT-M\"{u}ller and MT-PMCHWT formulations. \textit{Left}: Condition number as a function of the wavenumber $\kappa_0$, with the mesh size fixed at $h = 0.1\mathrm{m}$. \textit{Right}: Number of GMRES iterations required to reach a relative tolerance of $10^{-6}$ as a function of the mesh size $h$, with $\kappa_0 = 1\mathrm{m^{-1}}$. The condition number and GMRES iteration count of the MT-M\"{u}ller equation remain constant as $\kappa_0$ and $h$ tend to zero. In contrast, the MT-PMCHWT formulation exhibits rapidly deteriorating conditioning, with both metrics increasing significantly under mesh refinement and at low frequencies.}
    \label{fig:conditioning}
\end{figure*}

Fig.~\ref{fig:conditioning} presents the condition numbers and GMRES iteration counts for the mixed MT-M\"{u}ller discretization \eqref{eq:discretization} and the MT-PMCHWT formulation \eqref{eq:dis_PMCHWT} as functions of the wavenumber $\kappa_0$ and the mesh size $h$. The results clearly demonstrate that the MT-M\"{u}ller equation remains well-conditioned as $\kappa_0$ and $h$ tend to zero. In contrast, the MT-PMCHWT formulation exhibits rapidly increasing condition numbers and iteration counts under mesh refinement and in the low-frequency regime. The non-increasing condition number at extremely low frequencies further indicates that the proposed MT-M\"{u}ller formulation does not admit any non-trivial nullspaces, even on toroidal surfaces and in topologically complex geometries such as the dual-torus domain.

\subsection{Computational Time}

In the final experiments, we investigate the computational performance of the mixed MT-M\"{u}ller formulation \eqref{eq:discretization} in comparison with the standard MT-PMCHWT formulation \eqref{eq:dis_PMCHWT} and its CP variants. Specifically, we consider the following CP MT-PMCHWT formulation  \cite{LBG+2022}:
\begin{equation}
    \label{eq:CP_PMCHWT}
    \T \, \G^{-1} \Pbb \, \uv = \T \, \G^{-1} \wh{\uv}^{inc},
\end{equation}
where $\T$ and $\G$ are block-diagonal matrices whose diagonal entries are given by
\[
    \left[\T\right]_{kk} = 
    \begin{pmatrix}
        \zrb & \T_k \\
        -\T_k & \zrb
    \end{pmatrix}, \qq
    \left[\G\right]_{kk} = 
    \begin{pmatrix}
        \G_k & \zrb \\
        \zrb & \G_k
    \end{pmatrix},
\]
for $k = 1, 2, \ldots, N$. Here, the local Gram matrix $\G_k$ represents the bilinear form $\inprod{\cdot, I \cdot}_{\times, k}$ on $\text{RT}(\Gm_{h, k}) \times \text{BC}(\Gm_{h, k})$, while $\T_k$ is the matrix representation of $\langle\cdot, \wt{T}_{kk} \cdot \rangle_{\times, k}$ on $\text{BC}(\Gm_{h, k}) \times \text{BC}(\Gm_{h, k})$. The Yukawa-type single-layer operators $\wt{T}_{kk}$ are obtained from the operators $T^{(0)}_{kk}$ by replacing the wavenumber $\kappa_0$ with the purely imaginary value $-\iota \kappa_0$ \cite{LC2024}. 

\begin{figure*}[http]
    \centering
    \includegraphics[width=\linewidth]{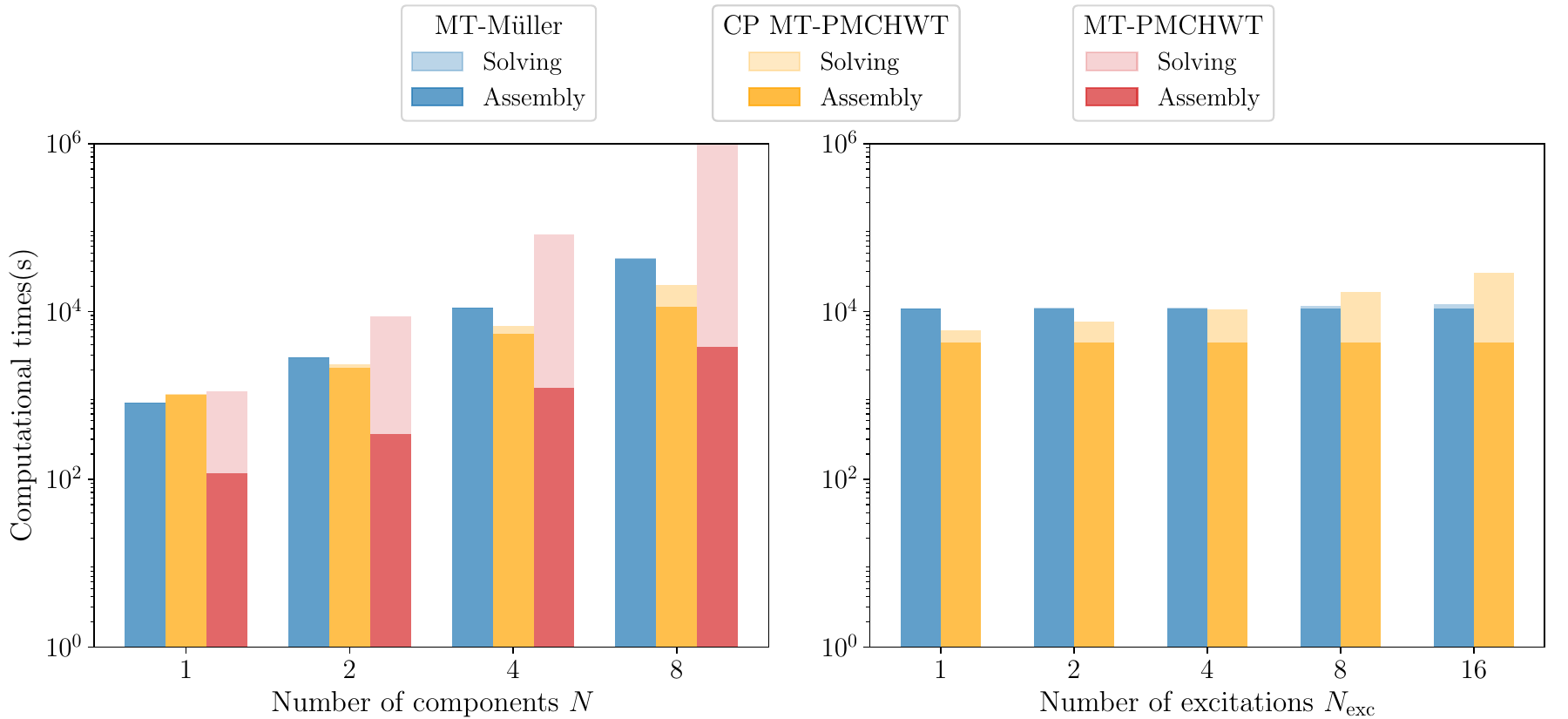}
    \caption{Assembly and solving times for different boundary integral formulations applied to electromagnetic scattering by a vertical stack of unit cubes, shown as functions of the number of cubes $N$ and the number of excitations $N_{\text{exc}}$. Each cube is discretized using a uniform triangulation with mesh size $h = 0.07\mathrm{m}$, resulting in $N_{e,k} = 4890$ edges (corresponding to $9780$ unknowns) per cube. \textit{Left}: varying numbers of cubes illuminated by a single plane wave. \textit{Right}: a fixed stack of four cubes excited by multiple plane waves.}
    \label{fig:cost}
\end{figure*}

Fig.~\ref{fig:cost} summarizes the measured assembly and solving times for the three formulations as the number of components $N$ and the number of excitations $N_\text{exc}$ increase. In both experiments, the incident fields are plane waves with background wavenumber $\kappa_0 = 0.01\mathrm{m}^{-1}$. In the first experiment, scatterers consisting of varying numbers of unit cubes are illuminated by a single plane wave propagating in the $\hat{\zb}$ direction and polarized along $\hat{\xb}$. In the second experiment, the scatterer is fixed as a stack of four cubes, while the $N_{\text{exc}}$ excitations are plane waves propagating in the $\hat{\zb}$ direction with polarization vectors
\[
    \hat{\pb}_m = \cos\paren{\dfrac{2\pi m}{N_{\text{exc}}}} \hat{\xb} + \sin \paren{\dfrac{2\pi m}{N_{\text{exc}}}} \hat{\yb}, 
\]
for all $m = 1, 2, \ldots, N_{\text{exc}}$, where $\hat{\yb}$ is the unit vector along the $y$-axis.

The MT-PMCHWT formulation exhibits minimal assembly costs but suffers from prohibitively large solution times, making preconditioning crucial. The CP MT-PMCHWT formulation \eqref{eq:CP_PMCHWT} significantly reduces the solution time at the expense of additional assembly work associated with the construction of the Calder\'{o}n preconditioner, resulting in a substantially improved overall performance. This overhead is already mitigated by employing a diagonal Calder\'{o}n preconditioner and would increase further if full Calder\'{o}n preconditioners were used instead.

In contrast, the mixed MT-M\"{u}ller formulation \eqref{eq:discretization} incurs higher assembly costs due to barycentric mesh refinement and the increased number of BIOs, but consistently achieves the shortest solution times. As a result, its total runtime is comparable to that of the CP MT-PMCHWT formulation. The short solving times are particularly advantageous in applications involving multiple right-hand sides. Moreover, in low-frequency regimes where the MT-PMCHWT formulations experience deteriorating conditioning, the mixed MT-M\"{u}ller formulation remains robust.

\section{Conclusion}
\label{sec:conclusion}

We have introduced a global multi-trace M\"{u}ller boundary integral formulation for electromagnetic scattering by composite dielectric objects, composed exclusively of second-kind operators. Employing a conforming mixed discretization based on RWG and BC basis functions, the method delivers accurate trace solutions and derived electromagnetic quantities while exhibiting favorable conditioning on dense meshes and in low-frequency regimes. From a computational perspective, the mixed MT-M\"{u}ller formulation achieves significantly reduced solution times. Although it entails increased assembly costs, its overall computational time remains comparable to that of Calder\'{o}n-preconditioned PMCHWT formulations in typical configurations. In scenarios involving multiple excitations, where solution efficiency dominates the total runtime, the MT-M\"{u}ller formulation offers a particularly compelling advantage, establishing it as a highly competitive alternative for large-scale and repeated-scattering simulations.

As already mentioned in the introduction, an important direction for future research concerns the extension of the proposed formulation to time-domain problems. Future work will therefore focus on a systematic investigation of TD-M\"{u}ller formulations and on extending the multi-trace framework to transient electromagnetic scattering.


%




\ifCLASSOPTIONcaptionsoff
  \newpage
\fi



\bibliographystyle{IEEEtran}
\bibliography{abrv_ref.bib}

%

\begin{IEEEbiography}[{\includegraphics[width=1in,height=1.25in,clip,keepaspectratio]{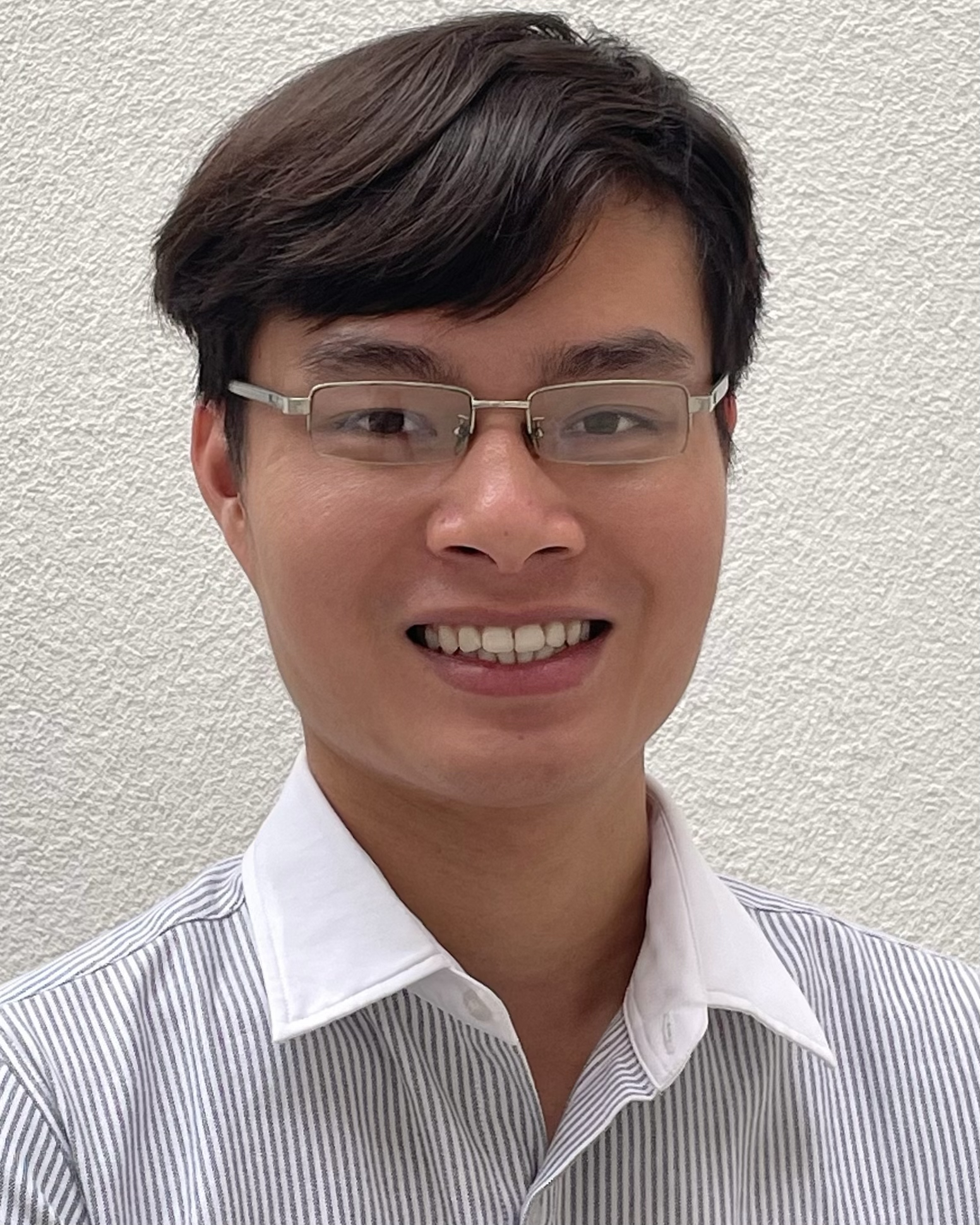}}]{Van Chien Le}(Member, IEEE)
received the M.Sc. degree in applied mathematics from Hanoi University of Science and Technology, Vietnam, in 2018, and the Ph.D. degree in mathematics from Ghent University, Belgium, in 2022.

Since 2022, he has been a Postdoctoral Researcher at the Department of Information Technology, Ghent University. His research interests include the numerical analysis of finite element and boundary element methods, with a particular focus on robust and stable frequency-domain and time-domain boundary integral equation solvers for electromagnetic scattering.

Dr. Le was awarded the IEEE Ulrich L. Rohde Innovative Conference Paper Award on Computational Techniques in Electromagnetics in 2023 and the URSI Young Scientist Award in 2026.
\end{IEEEbiography}

\begin{IEEEbiography}
[{\includegraphics[width=1in,height=1.25in,clip,keepaspectratio]{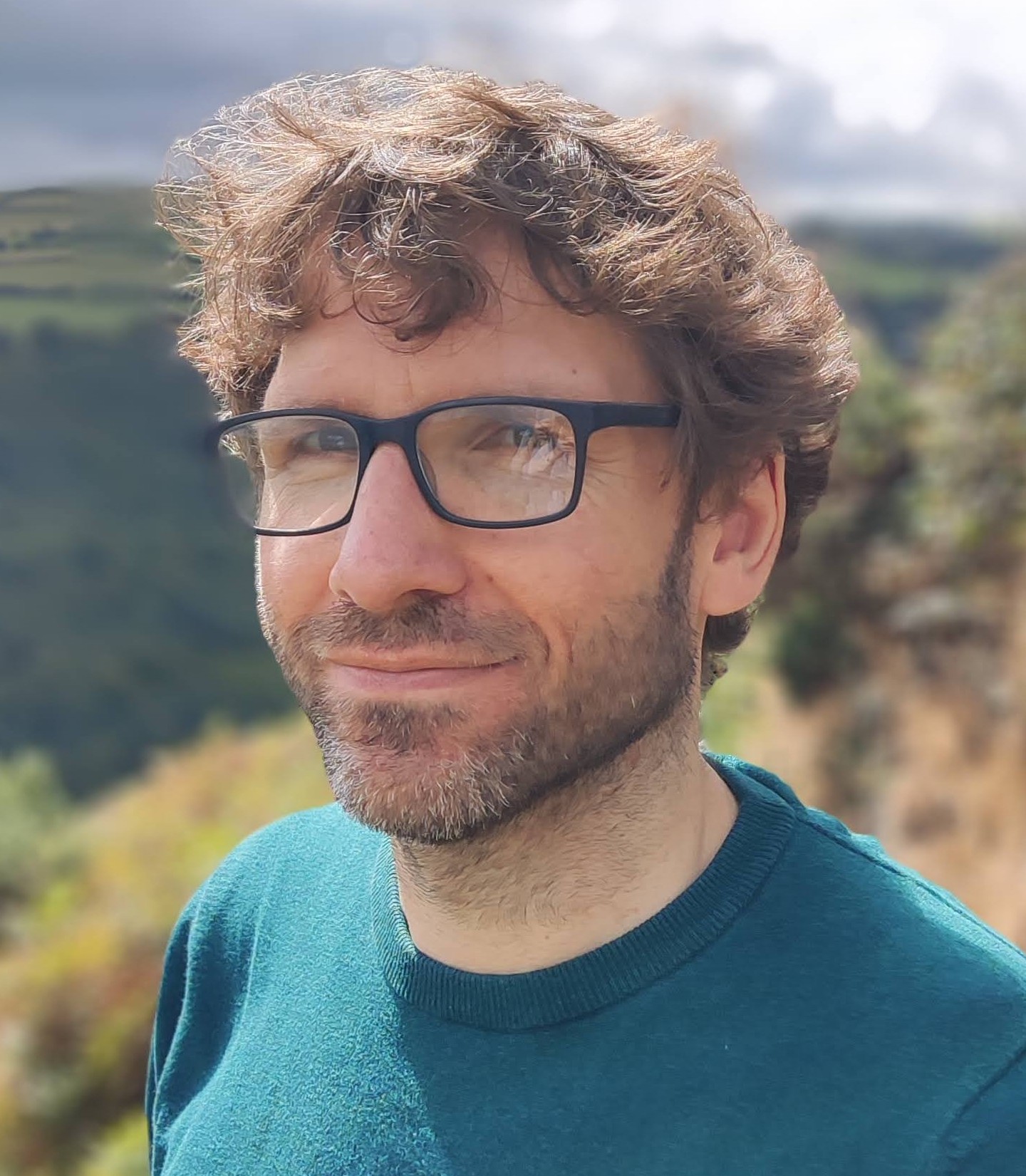}}]{Kristof Cools}(Member, IEEE) is full professor with the Department of Information Technology at Ghent University. Before joining Ghent University he has held positions at the TU Delft and the University of Nottingham. Kristof Cools has been awarded the ICEAA Young Scientist Award 2009, the URSI Young Scientist Award 2014, and the Rohde Most Innovative Conference Paper Award 2023. He has authored over 33 papers published in peer reviewed international journals and over 120 conference contributions, resulting in a Google Scholar recorded h-index of 19 and over 1800 citations. Dr. Cools has received funding from institutional, regional, industrial, national, and European funding bodies, including a Marie Skłodowska-Curie career integration grant and ERC consolidator/PoC project on boundary element methods for time-domain domain decomposition methods. 
\end{IEEEbiography}







\end{document}

%% file: lvc_style.tex
\newcommand{\R}{\mathbb R}
\newcommand{\C}{\mathbb C}

\newcommand{\G}{\mathbb G}
\newcommand{\M}{\mathbb M}

\newcommand{\T}{\mathbb T}
\newcommand{\Pbb}{\mathbb P}

\newcommand{\TT}{\mathcal T}

\newcommand{\KK}{\mathcal K}

\DeclareMathOperator{\diag}{diag}

\newcommand{\paren}[1]{\left({#1}\right)}

\newcommand{\abs}[1]{\left\vert{#1}\right\vert}

\newcommand{\inprod}[1]{\left\langle{#1}\right\rangle}

\newcommand{\wt}[1]{\widetilde{#1}}
\newcommand{\wh}[1]{\widehat{#1}}
\newcommand{\ovl}[1]{\overline{#1}}

\newcommand{\transpose}{\mathsf{T}}

\newcommand{\bm}[1]{\boldsymbol{#1}}

\newcommand{\curlt}{\textbf{\textup{curl}}}
\newcommand{\divt}{\textup{div}}
\newcommand{\gradt}{\textbf{grad}}

\newcommand{\eps}{\epsilon}

\newcommand{\pa}{\partial}

\newcommand{\Gm}{\Gamma}
\newcommand{\om}{\omega}
\newcommand{\Om}{\Omega}

\newcommand{\xb}{\bm{x}}
\newcommand{\yb}{\bm{y}}
\newcommand{\zb}{\bm{z}}

\newcommand{\hb}{\bm{h}}

\newcommand{\eb}{\bm{e}}

\newcommand{\ub}{\bm{u}}

\newcommand{\jb}{\bm{j}}

\newcommand{\mb}{\bm{m}}
\newcommand{\qb}{\bm{q}}

\newcommand{\pb}{\bm{p}}

\newcommand{\nv}{\textbf{n}}
\newcommand{\uv}{\textbf{u}}
\newcommand{\fv}{\textbf{f}}
\newcommand{\gv}{\textbf{g}}

\newcommand{\zrb}{\bm{0}}

\newcommand{\ds}{\,\mathrm{d}s}

\newcommand{\q}{\quad}
\newcommand{\qq}{\qquad}
\newcommand{\qqq}{\qquad\quad}
\newcommand{\qqqq}{\qquad\qquad}
\newcommand{\qqqqq}{\qquad\qquad\quad}